\documentclass[12pt,a4paper,equation]{article}
\usepackage[greek,american]{babel}
\usepackage[iso-8859-7]{inputenc}
\usepackage{amssymb,amsfonts}
\usepackage[dvips]{graphicx}
\usepackage{amsmath}
\usepackage{epsfig}
\usepackage{latexsym}
\usepackage{makeidx}
\usepackage{pst-math,pst-xkey}

\numberwithin{equation}{section}
\makeindex
\begin{document}

\date{}
\author{Vassilis G. Papanicolaou$^1$, and Kyriakie Vasilakopoulou$^2$\\
%EndAName
Department of Mathematics\\
National Technical University of Athens\\
Zografou Campus\\
157 80 Athens, GREECE\\
$^1$papanico@math.ntua.gr\\
$^2$vkiriaki@math.uoa.gr}

\title{Similarity Solutions of a Replicator Dynamics Equation Associated to a Continuum of Pure Strategies}
\maketitle
\begin{abstract}
We introduce a nonlinear degenerate parabolic equation containing a nonlocal term. The equation serves as a replicator dynamics model where the set
of strategies is a continuum. In our model the payoff operator (which is the continuous analog of the payoff matrix) is nonsymmetric and, also,
evolves with time. We are interested in solutions $u(t, x)$ of our equation which are positive and their integral (with respect to $x$)
over the whole space is $1$, for any $t > 0$. These solutions, being probability densities, can serve as time-evolving mixed
strategies of a player. We show that for our model there is an one-parameter family of self-similar such solutions $u(t, x)$, all approaching
the Dirac delta function $\delta(x)$ as $t \to 0^+$.
\end{abstract}


\textbf{2010 AMS Mathematics Classification.} 35C06; 91A22; 91B52.\\\\

\noindent\author{\underline{Corresponding author}: Vassilis G. Papanicolaou}

\underline{e-mail:} papanico@math.ntua.gr

\section{Introduction}

The replicator dynamics models are popular models in
evolutionary game theory. They have significant applications in
economics, population biology, as well as in other areas of science \cite{H-S}, \cite{I}, \cite{S}, \cite{T-J}.

Replicator dynamics have been studied extensively in the finite
dimensional case: Let $A = (a_{ij})$ be an $m \times m$ matrix. The typical
replicator dynamics equation is \cite{H-S}
\begin{equation}
u_t = \left[ Au - (u, Au) \right] u = (Au) u - (u, Au) u,
\label{A1}
\end{equation}
where the subscript $t$ in $u_t$ denotes derivative with respect to the time variable $t$
and $(Au) u$ is the vector whose $i$-th component is the product
of the $i$-th components of $(Au)$ and $u$ (i.e. the ``pointwise product" of two vectors).
The matrix $A$ is called the payoff matrix, while $ S = \{ 1,...,m \} $ is the
set of pure strategies (or options) and the vector
\begin{equation*}
u = \left( u_1(t),...,u_m(t) \right)^{\top},
\end{equation*}
is a probability (mass) function on $S$, meaning that
\begin{equation}
u_j(t) \geq 0, \quad \text{for } j = 1,...,m,
\qquad \text{and} \qquad
\sum_{j = 1}^m u_j(t) = 1.
\label{A1a}
\end{equation}

It is easy to see that if the conditions \eqref{A1a} are satisfied
for $t = 0$, then they are satisfied for all $t \geq 0$ (under the flow
(\ref{A1})). The vector $u$ represents the mixed strategy of one member of the population, i.e. one player, against the rest of the population.
The dependence of $u$ in $t$ allows the player to update her strategy, in order to increase her payoff.

Infinite dimensional versions of this evolutionary strategy models
have been proposed, e.g., in \cite{Bo} and \cite{Oe} (see also \cite{Oe1} and the
survey \cite{H-S}) in connection to certain economic and biological
applications. For instance, there are situations where (pure) strategies correspond to
geographical points and hence it is natural to model the set of
strategies by a continuum. However, the abstract form of the proposed equations does
not allow one to obtain much insight, for example on the form of solutions.

In order to make some progress in this direction, the recent works
\cite{K-P} and \cite{P-S} initiated the study of the case where
$S$ is the set $\Bbb{R}^d$, $d \geq 1$, and the payoff operator $A$
is the Laplacean operator $\Delta$. Then the evolution law (\ref{A1}) becomes
\begin{equation}
u_t = \left[ \Delta u - (u, \Delta u) \right] u,
\label{A2}
\end{equation}
where $(\cdot \, , \cdot)$ denotes the usual inner product of the
Hilbert space $L^2(\mathbb{R}^d)$ of the square-integrable
functions defined on $\mathbb{R}^d$.
References \cite{K-P} and \cite{P-S} deal only with the special problem of constructing
an one-parameter family of self-similar solutions for (\ref{A2}),
namely solutions $u$ of the form
\begin{equation*}
u(t,x) = t^{-\kappa} g\left(r t^{-\lambda}\right),
\qquad \text{where}\quad  r := |x| = \sqrt{x_1^2 + \cdots + x_d^2}.
\end{equation*}
A peculiar feature of these solutions is that all of them are probability densities on
$\mathbb{R}^d$, for all $t > 0$,
and approach the Dirac delta function $\delta(x)$ as $t \to 0^+$.

One criticism towards (\ref{A2}) is that the Laplacean operator $\Delta$ is a symmetric operator and, also, time-independent.
A payoff operator $A$ which is symmetric with respect to the inner product $(\cdot \, , \cdot)$
corresponds to the case of a {\em partnership game}, where interests of both players coincide (see, e.g., \cite{H-S}).
These are unrealistic features for a payoff operator in a replicator dynamics model.
For this reason, in the present work we consider a nonsymmetric and time-dependent payoff operator, namely
\begin{equation}
A u = A(t) \, u = \frac{\partial^2 u}{\partial x^2} + a t^{\gamma} x \frac{\partial u}{\partial x},
\label{A2a}
\end{equation}
where $\gamma$ is a specific constant (we will see later that $\gamma = -2/3$), while
\begin{equation}
a > 0
\label{A2aa}
\end{equation}
is an arbitrary but fixed constant. Then, our replicator dynamics model is described by the equation
\begin{equation}
u_t = \left[ Au - (u, Au) \right] u,
\qquad t > 0,
\quad
x \in \mathbb{R},
\label{A3}
\end{equation}
with $u = u(t, x)$, where the operator $A$ is given by \eqref{A2a}.

In order for (\ref{A3}) to be a replicator dynamics model, we need to make sure that if we start with an initial condition which is a probability
density function, namely
\begin{equation}
u(0, x) = f(x) \geq 0,
\quad x \in \mathbb{R},
\qquad \text{satisfying }\;
\int_{\mathbb{R}} f(x) dx = 1,
\label{A4}
\end{equation}
then the solution $u(t, x)$ will remain a probability density function (as a function of $x$) for all $t > 0$ (as long as it exists).
This can be justified as follows: Set
\begin{equation}
U(t) := \int_{-\infty}^{\infty} u(t, x) dx.
\label{J1}
\end{equation}
Then, integrating both sides of \eqref{A3} over $\mathbb{R}$ (with respect to $x$) gives
\begin{equation}
U'(t) = (u, Au) [1 - U(t)],
\label{J2}
\end{equation}
where we have assumed that the interchange of integration with respect to $x$ and differentiation with respect to $t$ is allowed.
Now, in view of \eqref{J2}, the fact that $U(0) = 1$ (which follows from \eqref{A4}) suggests that $U(t) \equiv 1$; in view of \eqref{J1}, 
this says that the integral of $u(t, x)$, with respect to $x$, on $\mathbb{R}$ is $1$ for every $t$.

Also, if $u(t, x)$ is a solution of \eqref{A3} which exists for all $t > 0$ and, as a function of $x$, it is integrable on $\mathbb{R}$ and
positive for small $t$, then, due to the nature of equation \eqref{A3} we have that $u(t, x)$ remains positive for all $t > 0$. We can, thus,
conclude that the set of probability densities on $\mathbb{R}$ is invariant under the flow \eqref{A3}.

%Also, if $u(t, x)$, as a function of $x$, attains its minimum at some $x_0 = x_0(t)$, and $u(t, x_0)$ approaches $0$, as $t \to t_0^{-}$, for some
%$t_0 > 0$ (so that $u(t_0, x_0) = 0$), then
%(\ref{A3}) together with the facts that $u_x(t_0, x_0) = 0$ and $u_{xx}(t_0, x_0) \geq 0$ give that $u_t(t_0, x_0) > 0$, contradicting the fact that %$u(t, x_0)$ moves down to $0$ as $t \to t_0^{-}$. This shows that if $u(0, x) > 0$, then $u(t, x) > 0$ for all $t > 0$
%(as long as it exists). 

\section{Special solutions}

\subsection{Self-Similar solutions}
We consider the equation (\ref{A3}), where $A$ is given by \eqref{A2a}. Let us assume that the solution $u(t, x)$ satisfies
\begin{equation}
u(t, \cdot \,) \in H^1(\mathbb{R})
\qquad \text{and} \qquad 
\lim_{x\rightarrow \pm \infty} x u(t,x)^2 = 0.
\label{5}
\end{equation}
By \eqref{A2a} we get
\begin{equation}\notag
(Au,u) = \int_{-\infty}^{\infty}(Au)udx
=\int_{-\infty}^{\infty}u_{xx}udx + a t^{\gamma} \int_{-\infty}^{\infty} x u_{x} u dx.
\notag
\end{equation}
Hence, in view of \eqref{5}, integration by parts yields
\begin{equation}\label{6}
 (Au,u) = - \int_{-\infty}^{\infty} u_x^2 dx - \frac {a}{2}t^{\gamma}\int_{-\infty}^{\infty}u^{2}dx,
\end{equation}
thus, (\ref{A3}) is equivalent to
\begin{equation}\label{7}
u_t = \left[u_{xx} + a t^{\gamma} xu_{x} + \int_{-\infty}^{\infty} u_x^2 dx + \frac{a}{2}t^{\gamma} \int_{-\infty}^{\infty}u^{2}dx \right] u.
\end{equation}
We will look for self-similar solutions of (\ref{A3}), namely solutions $u(t,x)$ of the form
\begin{equation}\label{8}
u(t,x) = t^{-\kappa}g(x t^{-\lambda}).
\end{equation}
We set $s=x t^{-\lambda}$ (hence $x=s t^{\lambda}$) so that $u(t,x)$ of \eqref{8} can be also written as
$u(t,x)=t^{-\kappa} g(s)$.
It follows that
\begin{equation}\label{9}
u_{x}(t,x) = t^{-(\kappa + \lambda)}g'(s)
\end{equation}
and
\begin{equation}\label{10}
u_{xx}(t,x) = t^{-(\kappa + 2\lambda)}g''(s)
\end{equation}
Also,
\begin{equation}\label{11}
u_{t}(t,x)=- \kappa t^{-(\kappa + 1)}g(s) - \lambda xt^{-\lambda}t^{-(\kappa + 1)}g'(s)
= -\kappa t^{-(\kappa + 1)}g(s) - \lambda st^{-(\kappa + 1)}g'(s).
\end{equation}
Then, \eqref{6} gives
\begin{equation}\label{12}
(Au,u) = - t^{-(2\kappa + \lambda)}\int_{-\infty}^{\infty} g'(s)^2 ds
- \frac{a}{2}t^{\gamma + \lambda - 2\kappa}\int_{-\infty}^{\infty} g(s)^2 ds.
\end{equation}
Setting
\begin{equation}\label{13}
K[g] := \int_{-\infty}^{\infty} g'(s)^2 ds
\qquad\text{and} \qquad
\Lambda[g] := \int_{-\infty}^{\infty} g(s)^2 ds,
\end{equation}
equation \eqref{12} becomes
\begin{equation}\label{14}
(Au,u)= - t^{-(2\kappa + \lambda)}K[g] - \frac{a}{2}t^{\gamma + \lambda - 2\kappa}\Lambda[g]
\end{equation}
Substituting \eqref{9}, \eqref{10}, \eqref{11}, \eqref{14} in \eqref{7}, we have
\begin{align}
-\kappa g(s) - \lambda s g'(s) = & \,t^{1 - \kappa - 2 \lambda}g''(s)g(s) + a s t^{1 + \gamma-\kappa}g'(s)g(s) + t^{1 - 2\kappa - \lambda} K[g]g(s)
\nonumber
\\
&+ \frac{a}{2}t^{1 + \gamma + \lambda - 2\kappa}\Lambda[g]g(s).
\label{15}
\end{align}
The only way that the above is a meaningful equation is that it does not contain $t$, which means that
\begin{equation}
1 - \kappa - 2 \lambda = 0, \quad 1 + \gamma-\kappa = 0,\quad 1 - 2\kappa - \lambda = 0,\quad 1 + \gamma + \lambda - 2\kappa = 0.
\end{equation}
This gives
\begin{equation}\label{17}
\gamma = - \frac{2}{3},
\quad
\kappa = \frac{1}{3},
\quad
\lambda=\frac{1}{3}.
\end{equation}
Finally, we notice that, under \eqref{17}, \eqref{8} gives
\begin{equation}\notag
\int_{-\infty}^{\infty}u(t,x)dx=\int_{-\infty}^{\infty}t^{-\kappa}g(x t^{-\lambda})dx=\int_{-\infty}^{\infty}t^{-1/3}g(x t^{-1/3})dx=\int_{-\infty}^{\infty}g(s)ds.
\end{equation}
which is independent of t. Thus, if we set
\begin{equation*}
\int_{-\infty}^{\infty}g(s)ds=1
%\label{A5}
\end{equation*}
then
\begin{equation}\notag
\int_{-\infty}^{\infty}u(t,x)dx=1, \text{ for all } t\geq0.
\end{equation}
The following lemma summarizes what we have done so far.

\smallskip

\textbf{Lemma 2.1.} If
\begin{equation}\label{18}
u(t,x) = t^{-\kappa} g(x t^{-\lambda})
\end{equation}
is a probability density in $x$ and satisfies \eqref{7}, then we must have
\begin{equation}\label{19}
\gamma = - \frac{2}{3},
\quad
\kappa = \frac{1}{3},
\quad
\lambda=\frac{1}{3},
\end{equation}
\begin{equation}\label{20}
g(s) \geq 0,
\quad s \in \mathbb{R},
\qquad
\int_{-\infty}^{\infty}g(s)ds = 1
\end{equation}
and
\begin{equation}\label{21}
g''(s)g(s) + a s g'(s)g(s) + K[g]g(s) + \frac{a}{2}\Lambda[g]g(s) + \frac{1}{3}g(s) + \frac{1}{3}sg'(s) = 0,
\end{equation}
where
\begin{equation}\label{22}
K[g] = \int_{-\infty}^{\infty} g'(s)^2 ds
\end{equation}
and
\begin{equation}\label{23}
\Lambda[g] = \int_{-\infty}^{\infty} g(s)^2 ds.
\end{equation}
Conversely, if \eqref{19}--\eqref{23} hold, then $u(t, x)$ given by \eqref{18} is a probability density in $x$ and satisfies \eqref{7}.

\smallskip

\textbf{Remark.} In view of \eqref{A2a}, the fact that $\gamma = -2/3$ tells us that in the long run and as long as $x$ stays bounded,
the payoff operator $A(t)$ of our model approaches the symmetric operator $\partial^2 / \partial x^2$.

\smallskip

Next, we need to show that there exist function(s) $g(s)$ satisfying \eqref{20} and \eqref{21}.

\subsection{The auxiliary problem}

Consider the problem
\begin{equation}\label{24}
q''(s)q(s) + a s q'(s)q(s) + \mu q(s) + \frac{1}{3}sq'(s) = 0
\end{equation}
\begin{equation}\label{25}
q(0) = A > 0,
\qquad
q'(0)=0
\end{equation}
where $\mu$ is a real parameter satisfying
\begin{equation}\label{26}
\mu >\frac{1}{3}.
\end{equation}
Equation \eqref{24} can be written in the form
\begin{equation}\label{27}
q''(s) + \left[\frac{1}{3q(s)} + a\right] s q'(s) + \mu = 0
\end{equation}
as long as $q(s)\neq 0$. Since $q(0) = A > 0$, the standard existence and uniqueness theorems for ordinary differential equation
imply that there is a $\delta > 0$ such that \eqref{24}-\eqref{25} has a unique solution $q(s)$ for $s \in (-\delta, \delta)$.
In fact, due to the invariance of \eqref{24} under the transformation $s \mapsto - s$ and the fact that $q'(0) = 0$, we must have
\begin{equation}\notag
q(-s) = q(s),
\qquad
s \in (-\delta, \delta).
\end{equation}

\smallskip

\textbf{Lemma 2.2.}
The solution q(s) of \eqref{24}-\eqref{25} exists for all $s \in \mathbb{R}$ and it is a strictly positive (even) function which is decreasing on
$(0, \infty)$. Also,
\begin{equation}\label{28}
\lim_{s\rightarrow \infty}q(s) = \lim_{s\rightarrow \infty}q'(s) = 0,
\end{equation}
\begin{equation}\label{29}
\int_{-\infty}^{\infty} q'(s)^2 ds < \infty,
\end{equation}
and
\begin{equation}\label{30}
\int_{-\infty}^{\infty} q(s)^2 ds <\infty.
\end{equation}
Furthermore, the following equality holds
\begin{equation}\label{31}
\left(\mu - \frac{1}{3}\right)\int_{-\infty}^{\infty} q(s)ds
= \int_{-\infty}^{\infty} q'(s)^2 ds + \frac{a}{2}\int_{-\infty}^{\infty} q(s)^2 ds.
\end{equation}

\smallskip

\textit{Proof}.
Since $q$ is an even function, it is enough to show that $q(s)$ exists for all $s \in [0, \infty)$.
If this is not true, then either
(i) (due the denominator $q(s)$ appearing in \eqref{27}) there must be an $s_{1} \in (0 , \infty)$ such that
$q(s_{1})=0$, while $q(s) > 0$ for all $s \in [0 , s_{1})$, or
(ii) by a well-known theorem in the theory of ordinary differential equations \cite{C-L} there must exist some $b > 0$
such that
\begin{equation*}
\lim_{s \rightarrow b^{-}}\left[ \left|q'(s)\right| + \left| q(s)\right|\right]= \infty.
%\label{32}
\end{equation*}

Let us first exclude the case (i).
Suppose that there is an $s_{1} > 0$ such that $q(s_{1}) = 0$, while $q(s) > 0$ for all $s \in [0, s_{1})$.
Then, $q'(s)$ is negative in $(0, s_{1})$.
If this were not true, then there should exist a $s_{2} \in (0, s_{1})$ such that $q'(s_{2}) = 0$, while $q'(s) < 0$ for all $s\in (0, s_{2})$.
This would imply that $q''(s_2) \geq 0$. However, by \eqref{27}
\begin{equation}\notag
q''(s_{2})=-\mu<0,
\end{equation}
a contradiction.

Now, if we integrate \eqref{27} from $0$ to $s \in (0, s_{1})$  and use the fact that $q'(0)=0$, we get
\begin{equation}\notag
\int_{0}^{s}\left[ q''(\xi) + \frac{\xi q'(\xi)}{3q(\xi)} + a \xi q'(\xi) + \mu \right]d\xi
= q'(s) + \frac{1}{3}\int_{0}^{s} \xi \left[\ln q(\xi)\right]'d\xi +a sq(s) - a \int_{0}^{s}q(\xi)d\xi + \mu s = 0,
\end{equation}
or
\begin{equation}\label{33}
q'(s) = a \int_{0}^{s}q(\xi)d\xi - a sq(s)  - \mu s + \frac{1}{3}\int_{0}^{s}\ln q(\xi)d\xi - \frac{1}{3} s \ln q(s).
\end{equation}
Since $q(s) > 0$, while $q'(s) < 0$ for $s \in (0, s_{1})$, $q(s)$ is decreasing in $[0, s_{1})$ and, consequently, $\ln q(s)$ is decreasing in
$[0, s_{1})$. Hence, the function 
\begin{equation}
f(s) := - \frac{1}{3} \ln q(s)
\end{equation}\label{333a}
is increasing in $(0, s_{1})$ and
\begin{equation}\notag
\lim_{s \rightarrow s_{1}^{-}}f(s) = -\frac{1}{3} \lim_{s \rightarrow s_{1}^{-}} [\ln q(s)] = \infty.
\end{equation}
Then, it is not hard to show (see, e.g., Proposition A.1 of the Appendix of \cite{K-P}) that
\begin{equation}\notag
\lim_{s \rightarrow s_{1}^{-}} \left(sf(s) - \int_{0}^{s} f(\xi)d\xi \right)  = \infty,
\end{equation}
i.e. (recall \eqref{333a})
\begin{equation}\notag
\lim_{s \rightarrow s_{1}^{-}} \left( - \frac{1}{3}s \ln q(s) + \frac{1}{3}\int_{0}^{s} \ln q(\xi) d\xi \right)  = \infty.
\end{equation}
Hence, \eqref{33} gives
\begin{equation}\notag
\lim_{s \rightarrow s_{1}^{-}}q'(s) 
= \lim_{s \rightarrow s_{1}^{-}} \left[a \int_{0}^{s}q(\xi)d\xi - a sq(s)  
- \mu s + \frac{1}{3}\int_{0}^{s}\ln q(\xi)d\xi - \frac{1}{3} s \ln q(s)\right]
= \infty,
\end{equation}
which is impossible, since, as we have seen, $q'$ stays negative in $(0, s_{1})$.
Hence such an $s_{1}$ cannot exist, i.e. $q$ never vanishes and consequently, $q'$ also never vanishes.
In particular, $q(s) > 0$, $q'(s) < 0$ (hence, $q$ is decreasing), and, therefore, $0 < q(s) < q(0) = A$, for all $s > 0$
for which $q(s)$ and $q'(s)$ exist.

Now suppose that there is an $b>0$ such that
\begin{equation}\label{34}
\lim_{s \rightarrow b^{-}}\left[ \left|q'(s)\right| + \left|q(s)\right| \right] = \infty.
\end{equation}
By the previous discussion, the only way for \eqref{34} to happen is
\begin{equation}\notag
\lim_{s \rightarrow b^{-}}q'(s) = - \infty.
\end{equation}
Then,
\begin{equation}\notag
\liminf_{s \rightarrow b^{-}}q''(s) = -\infty,
\end{equation}
which contradicts \eqref{27}.
Thus $q'$ remains finite and strictly negative on $(0, \infty)$ while $q$ is strictly positive and strictly decreasing
on $(0, \infty)$. Due to the evenness of $q$, we must have also $q(s) > 0$ for all $s < 0$. Hence, $q(s) > 0$ for all $s \in \mathbb{R}$.

From the previous discussion it follows that
\begin{equation}\notag
\lim_{s \rightarrow \infty} q(s) = L,
\end{equation}
namely
\begin{equation}\label{35}
q(s) = L + o(1)
\qquad
\text{as }\; s \rightarrow \infty,
\end{equation}
where $L \in [0, A)$. To continue, let us suppose $L > 0$.
Then, the above formula implies that, as $s \rightarrow \infty$,
\begin{equation}\label{36}
\ln q(s) = \ln \left(L + o(1)\right) = \ln L \left(1 + o(1)\right) = \ln L + o(1).
\end{equation}
Using \eqref{35} and \eqref{36} in \eqref{33}, we obtain
\begin{equation}\notag
q'(s) = a \int_{0}^{s} \left[ \ln L + o(1)\right] d\xi - a s \left[ \ln L + o(1)\right] 
-\mu s + \frac{1}{3}\int_{0}^{s} \left[\ln L + o(1)\right] d\xi - \frac{1}{3} s \left[\ln L + o(1)\right],
\end{equation}
which implies
\begin{equation}\notag
q'(s) = a s \ln L + o(s) - a s \ln L + o(s) - \mu s + \frac{1}{3} s \ln L + o(s) - \frac{1}{3} s \ln L + o(s),
\end{equation}
i.e.
\begin{equation}\notag
q'(s) = - \mu s + o(s)
\qquad
\text{as }\; s \rightarrow \infty,
\end{equation}
which contradicts \eqref{35}. Therefore $L = 0$, i.e.
\begin{equation}
\lim_{s \rightarrow \infty} q(s) = 0.
\end{equation}
We continue by noticing that
\begin{equation}\label{37}
\int_{0}^{\infty}q'(s)ds = \lim_{s \rightarrow \infty} q(s) - q(0) = - A,
\end{equation}
hence $q' \in L_{1}(\mathbb{R})$ (since $q'$ is odd and negative).
Suppose
\begin{equation}\label{38}
\liminf_{s \rightarrow \infty}q'(s) < 0.
\end{equation}
From \eqref{37}, there is a sequence $s_{n} \rightarrow \infty$ such that
$q'$ attains a local minimum at $s_{n}$ and
\begin{equation}\label{39}
\lim_{n \rightarrow \infty}q'(s_{n}) = - \delta, \text{ for some } \delta >0.
\end{equation}
But, since $q'(s_{n})$ is a local minimum we must have $q''(s_{n}) = 0$, hence \eqref{27} gives
\begin{equation}\notag
\left[\frac{1}{3q(s_{n})} + a\right] s_{n} q'(s_{n}) = - \mu
\qquad \text{or} \qquad 
q'(s_{n}) = - \frac{3 \mu q(s_{n})}{\left[1 + 3 a q(s_{n})\right] s_{n}},
\end{equation}
thus
\begin{equation*}
\lim_{n \rightarrow \infty}q'(s_{n}) = 0,
\end{equation*}
contradicting \eqref{39} and hence \eqref{38}.
We have, thus, established that
\begin{equation}\label{40}
\lim_{s \rightarrow \infty}q'(s) = 0.
\end{equation}
This, together with the fact that $q'$ is odd and integrable, implies $q' \in L_{2}(\mathbb{R})$, i.e.
\begin{equation}\label{42}
\int_{- \infty}^{\infty} q'(s)^2 ds < \infty.
\end{equation}
Finally, \eqref{27} implies
\begin{equation}\notag
\int_{0}^{s}q(\xi)q''(\xi)d\xi + \frac{1}{3}\int_{0}^{s} \xi q'(\xi)d\xi + a \int_{0}^{s} \xi q(\xi)q'(\xi)d\xi + \mu\int_{0}^{s}q(\xi)d\xi =0.
\end{equation}
By integrating by parts the first two terms above and using the fact that $q'(0) = 0$ we have
\begin{equation}\notag
q(s)q'(s) - \int_{0}^{s} q'(\xi)^2 d\xi + \frac{1}{3}sq(s) + \frac{a}{2}sq^{2}(s) - \frac{a}{2}\int_{0}^{s}q(\xi)^2 d\xi 
+ \left(\mu - \frac{1}{3}\right)\int_{0}^{s}q(\xi)d\xi = 0.
\end{equation}
Letting $s \rightarrow \infty$, the above equation implies
\begin{equation}\label{43}
\left(\mu - \frac{1}{3}\right)\int_{0}^{\infty} q(\xi)d \xi \, + \, \lim_{s \rightarrow \infty} \left[\frac{1}{3}sq(s) + \frac{a}{2}s q(s)^2\right]
= \int_{0}^{\infty} q'(\xi)^2 d\xi \, + \,\frac{a}{2}\int_{0}^{\infty}q(\xi)^2 d\xi
\end{equation}
Since $a > 0$ and $q(s) > 0$, \eqref{43} gives
\begin{equation}\label{44}
\left(\mu - \frac{1}{3}\right)\int_{0}^{\infty}q(\xi) d\xi 
\leq \int_{0}^{\infty} q'(\xi)^2 d\xi + \frac{a}{2} \int_{0}^{\infty} q(\xi)^2 d\xi.
\end{equation}
\begin{enumerate}
\item If we suppose that
\begin{equation}\notag
\int_{0}^{\infty} q(\xi) d\xi= \infty,
\end{equation}
then \eqref{44}, due the \eqref{42}, implies that
\begin{equation}\notag
\int_{0}^{\infty}q(\xi)^2 d\xi = \infty.
\end{equation}
setting
\begin{equation}\notag
M := \int_{0}^{\infty} q'(\xi)^2 d\xi < \infty,
\end{equation}
formula \eqref{44} implies
\begin{equation}\notag
\left(\mu - \frac{1}{3}\right)\int_{0}^{\infty} q(\xi)d\xi \leq M + \frac{a}{2}\int_{0}^{\infty}q(\xi)^2 d\xi.
\end{equation}
Furthermore, $q(\xi)>0$, for all $\xi \in (0, \infty)$, thus $\int_{0}^{\infty}q(\xi)d\xi > 0$ and then from the
above inequality we have
\begin{equation}\label{45}
\mu - \frac{1}{3}\leq \frac{M + \frac{a}{2}\int_{0}^{\infty} q(\xi)^2 d\xi}{\int_{0}^{\infty}q(\xi)d\xi}.
\end{equation}
Then
\begin{equation}\notag
\lim_{s\rightarrow \infty}\frac{M +\frac{a}{2}\int_{0}^{\infty}q(\xi)^2 d\xi}{\int_{0}^{\infty}q(\xi)d\xi}= \lim_{s\rightarrow \infty}\frac{\frac{a}{2}q(s)^2}{q(s)}=\frac{a}{2}\lim_{s\rightarrow \infty}q(s)=0
\end{equation}
But, then, from \eqref{45} we have $\mu \leq 1/3$, which contradicts the fact that $\mu > 1/3$.

Consequently,
\begin{equation}\label{46}
\int_{0}^{\infty}q(\xi)d\xi< \infty.
\end{equation}

\item The function $q$ is strictly positive and strictly decreasing on $(0, \infty)$ with $0<q(s)<A$, for all $s>0$.
Since \begin{equation}\notag
\lim_{s \rightarrow \infty}q(s)=0,
\end{equation}
there is a $s_{0} > 0$ such that $0 < q(s)^2 < q(s)$ for all $s \geq s_{0}$. Hence
\begin{equation}\notag
\lim_{s \rightarrow \infty} q(s)^2 = 0
\qquad \text{and} \qquad
0 < \int_{s_{0}}^{\infty}q(s)^2 ds < \int_{s_{0}}^{\infty} q(s) ds.
\end{equation}
Thus, from \eqref{46} it follows that
\begin{equation}\label{47}
\int_{0}^{\infty}q(s)^2 ds < \infty.
\end {equation}
\end{enumerate}
From \eqref{43}, \eqref{42}, \eqref{46}, and \eqref{47} we have
\begin{equation}\notag
\lim_{s \rightarrow \infty} \left[\frac{1}{3}sq(s) +\frac{a}{2} s q(s)^2\right] = \int_{0}^{\infty} q'(\xi)^2 d\xi 
+ \frac{a}{2}\int_{0}^{\infty}q(\xi)^2 d\xi
- \left(\mu - \frac{1}{3}\right)\int_{0}^{\infty}q(\xi)d\xi < \infty.
\end{equation}
Thus,
\begin{equation}\notag
\lim_{s \rightarrow \infty} \left[\frac{1}{3}sq(s) +\frac{a}{2} sq^{2}(s)\right] = L' \in \mathbb{R}.
\end{equation}
If $L' \neq 0$, then the above limit tells us that $T(s) := (1/3) q(s) + (a/2) q(s)^2$ is asymptotic to $L'/s$, contradicting the fact
that $T(s)$ is integrable. Therefore,
\begin{equation}\notag
L'=0.
\end{equation}
Then, \eqref{43} gives
\begin{equation}\notag
\left(\mu - \frac{1}{3}\right)\int_{0}^{\infty}q(s)ds = \int_{0}^{\infty} q'(s)^2 ds + \frac{a}{2}\int_{0}^{\infty} q(s)^2 ds,
\end{equation}
from which \eqref{31} follows immediately. The proof of this key lemma is now complete.
\hfill $\blacksquare$

\subsection{The construction of the self-similar solutions}

\textbf{Lemma 2.3.}
Let $q(s)$ be the solution of the problem \eqref{24}-\eqref{25}. Then
\begin{equation}\label{48}
\left\| q' \right\|_{\infty} \leq \mu \sqrt{\frac{3A}{1 + 3 a A}} \, ,
\end{equation}
where $\left\| \cdot \right\|_{\infty}$ denotes the sup-norm, as usual. Also
\begin{equation}\label{49}
\int_{0}^{\infty}q(s)ds \geq \frac{A^{3/2}\sqrt{1 + 3 a A}}{2\sqrt{3} \, \mu}
\end{equation}
and
\begin{equation}\label{50}
\int_{0}^{\infty}q(s)^2 ds \geq \frac{A^{5/2}\sqrt{1 + 3 a A}}{3\sqrt{3} \, \mu}.
\end{equation}

\smallskip

\textit{Proof}.
The function $q'$ is odd, hence
\begin{equation}\notag
\left\| q' \right\|_{\infty} = \sup \left\{\left|q'(s)\right| : s \geq 0 \right\}.
\end{equation}
Since $q'(s)<0$ in $(0, \infty)$ with
\begin{equation}\notag
q'(0)=0=\lim_{s\rightarrow \infty} q'(s),
\end{equation}
it follows that $q'$ attains its absolute minimum at some $s_m \in (0, \infty)$, and hence
\begin{equation}\notag
\left\| q' \right\|_{\infty} = \sup \left\{- q'(s) : s \geq 0 \right\}= - q'(s_m)=\left| q'(s_m)\right|.
\end{equation}
Also, $q''(s_m)=0$, thus \eqref{27} implies
\begin{equation}\notag
q'(s_m)= - \frac{\mu}{\left[\frac{1}{3q(s_m)} + a\right] s_m}
\end{equation}
therefore
\begin{equation}\notag
\left\| q' \right\|_{\infty}= - q'(s_m) = \frac{\mu}{\left[\frac{1}{3q(s_m)} + a\right] s_m}.
\end{equation}
But $q(s)$ is decreasing in $[0, \infty)$, while $q(0)=A$ and $s_m \in (0, \infty)$, hence
\begin{equation}\notag
\frac{\mu}{[\frac{1}{3q(s_m)} + a]s_m} \leq \frac{\mu}{[\frac{1}{3A} + a]s_{m}} < \frac{\mu}{a s_{m}}
\end{equation}
and then
\begin{equation}\label{51}
\left\| q' \right\|_{\infty}\leq \frac{3 \mu A}{(1 + 3 a A)s_{m}}< \frac{\mu}{a s_{\mu}}.
\end{equation}
Also, by \eqref{27}
\begin{equation}\notag
q''(s) + \mu = - \left[\frac{1}{3q(s)}+ a \right]sq'(s)\geq 0
\qquad \text{for all }\; s \geq 0, \ \text{ while }\; q'(0)=0.
\end{equation}
Thus, we must have
\begin{equation}\notag
q'(s)\geq - \mu s
\qquad \text{for all }\; s\geq 0,
\end{equation}
in particular
\begin{equation}\label{52}
\left\| q' \right\|_{\infty}= - q'(s_m) \leq \mu s_{m}.
\end{equation}
By combining \eqref{51} and \eqref{52} we obtain
\begin{equation}\notag
\left\| q' \right\|_{\infty} \leq \min \left\{ \frac{3 \mu A}{(1 + 3 a A)s_{m}} \; , \; \mu s_{m}\right\}.
\end{equation}
But, no matter what $s_{m}$ is, the quantity $\min\left\{ 3 \mu A (1 + 3 a A)^{-1} s_m^{-1} , \mu s_{m}\right\}$
(since the first term is decreasing in $s_{m}$ while the second is increasing) is always at most $M'$, where
\begin{equation}\notag
M' := \mu s^{\ast} = \frac{3 \mu A}{(1+3 a A) s^{\ast}}.
\end{equation}
Then
\begin{equation}\notag
s^{\ast} = \sqrt{\frac{3A}{1+3 a A}}
\qquad \text{and} \qquad
M' = \mu \sqrt{\frac{3A}{1+3 a A}}.
\end{equation}
Thus,
\begin{equation}\notag
\left\| q' \right\|_{\infty} \leq M',
\end{equation}
which is \eqref{48}. Furthermore,
\begin{equation}\notag
\left\| q' \right\|_{\infty} \geq - q'(s)
\qquad \text{for all }\; s \geq 0,
\end{equation}
hence
\begin{equation}\notag
q(s) \geq q(0) - s \left\| q' \right\|_{\infty} = A - s \left\| q' \right\|_{\infty}
\qquad \text{for all }\; s \geq 0.
\end{equation}
Then, by \eqref{48} we have
\begin{equation}\label{53}
q(s) \geq A - s \mu \sqrt{\frac{3A}{1 + 3 a A}}
\qquad \text{for all }\; s \geq 0,
\end{equation}
in particular for
\begin{equation}\notag
0 \leq s \leq \frac{\sqrt {\left(1 + 3 a A \right) A}}{\mu \sqrt{3}},
\end{equation}
since $q(s) > 0$ for all $s \geq 0$. Then (see \eqref{53}),
\begin{equation}\notag
\int_{0}^{\infty}q(s)ds \geq \int_{0}^{\frac{\sqrt {\left(1 + 3 a A \right) A}}{\mu \sqrt{3}}}q(s)ds
 \geq \int_{0}^{\frac{\sqrt {\left(1 + 3 a A \right) A}}{\mu \sqrt{3}}} \left( A - s \mu \sqrt{\frac{3A}{1 + 3 a A}} \right)ds
 = \frac{A^{3/2}\sqrt{1 + 3 a A}}{2\sqrt{3} \, \mu},
\end{equation}
which is \eqref{49}.

Finally, from \eqref{53} we also have
\begin{equation}\notag
\int_{0}^{\infty} q(s)^2 ds \geq \int_{0}^{\frac{\sqrt {\left(1 + 3 a A \right) A}}{\mu \sqrt{3}}}q(s)^2 ds
 \geq \int_{0}^{\frac{\sqrt {\left(1 + 3 a A \right) A}}{\mu \sqrt{3}}} \left( A - s \mu \sqrt{\frac{3A}{1 + 3 a A}} \right)^{2}ds =
 \frac{A^{5/2}\sqrt{1 + 3 a A}}{3\sqrt{3} \, \mu},
\end{equation}
which is \eqref{50}.
\hfill $\blacksquare$

\smallskip

\textbf{Corollary 2.1.}
If $q(s)$ satisfies \eqref{24}-\eqref{25}, then
\begin{equation}\label{54}
\lim_{A \rightarrow \infty} \int_{- \infty}^{\infty}q(s)ds = \infty,
\end{equation}
\begin{equation}\label{55}
\lim_{A \rightarrow \infty} \int_{- \infty}^{\infty} q(s)^2 ds = \infty,
\end{equation}
and
\begin{equation}\label{56}
\lim_{A \rightarrow 0^+} \int_{- \infty}^{\infty} q'(s)^2 ds = 0.
\end{equation}

\smallskip

\textit{Proof}.
By \eqref{49} and the evenness of $q(s)$ we have
\begin{equation}\label{57}
\int_{\infty}^{\infty}q(s)ds = 2 \int_{0}^{\infty}q(s)ds \geq \frac{2 A^{3/2}\sqrt{1 + 3 a A}}{3\sqrt{3} \, \mu}
\end{equation}
and since
\begin{equation}\notag
\lim_{A \rightarrow \infty} \frac{2 A^{3/2}\sqrt{1 + 3 a A}}{3\sqrt{3} \, \mu} = \infty,
\end{equation}
we get that \eqref{57} implies that
\begin{equation}\notag
\lim_{A \rightarrow \infty} \int_{- \infty}^{\infty}q(s)ds = \infty.
\end{equation}
The function $q$ is even, and hence $q^{2}$ is even too.
Furthermore, \eqref{50} implies
\begin{equation}\label{58}
\int_{\infty}^{\infty} q(s)^2 ds = 2 \int_{0}^{\infty} q(s)^2 ds \geq \frac{2 A^{5/2}\sqrt{1 + 3 a A}}{3\sqrt{3} \, \mu},
\end{equation}
hence, from \eqref{58} we have
\begin{equation}\notag
\lim_{A \rightarrow \infty} \int_{- \infty}^{\infty} q(s)^2 ds = \infty.
\end{equation}
Recall that $- q'(s) > 0$ (and $- q'(s) \leq \left\|q'\right\|_{\infty}$) for all $s \in (0, \infty)$. Thus, by \eqref{28} we get
\begin{equation}\notag
0 \leq \int_{0}^{\infty} q'(s)^2 ds \leq - \left\|q'\right\|_{\infty} \int_{0}^{\infty} q'(s)ds =
A \left\|q'\right\|_{\infty}
\end{equation}
and, consequently, by using \eqref{48} we have
\begin{equation}\label{59}
0 \leq \int_{0}^{\infty} q'(s)^2 ds \leq A \mu \sqrt{\frac{3A}{1 + 3 a A}} = \sqrt{\frac{3A^{3} \mu^{2}}{1 + 3 a A}}.
\end{equation}
Finally, since $q'(s)$ is odd and hence $q'(s)^{2}$ is even, by using \eqref{59} we get
\begin{equation}\label{60}
0 \leq \int_{- \infty}^{\infty} q'(s)^2 ds = 2 \int_{0}^{\infty} q'(s)^2 ds \leq 2\sqrt{\frac{3A^{3} \mu^{2}}{1 + 3 a A}},
\end{equation}
which implies \eqref{56}.
\hfill $\blacksquare$

\smallskip

\textbf{Lemma 2.4.}
If $q(s)$ is the solution of \eqref{24}-\eqref{25}, then
\begin{equation}\label{61}
\frac{q(1) e^{3 a q(1)}}{s^{3 \mu} \exp \left[3 \mu + \frac{3 \mu}{s} \sqrt{\frac{3 A}{1 + 3 a A}}\right]} <q(s) e^{3 a q(s)} \leq \frac{A e^{3 A (1 + a)}}{s^{3 \mu} \exp \left[3 \mu  - \frac{3 \mu}{s} \sqrt{\frac{3 A}{1 + 3 a A}}\right]}
\end{equation}
for all $s \geq 1$.

\smallskip

\textit{Proof}.
We consider the function
\begin{equation*}
F(s) := - \frac{1}{3} \int_{0}^{s}\ln q(\xi)d\xi,
\qquad
s \in \mathbb{R}.
\end{equation*}
Since $q(s)$ is decreasing in $(0, \infty)$, for $0 \leq \xi \leq 1$ we have
\begin{equation}\label{62}
 q(1) \leq q(\xi) \leq q(0) = A
\end{equation}
and
\begin{equation}\label{63}
\ln q(1) \leq -3F(1) \leq \ln A.
\end{equation}
Furthermore, \eqref{62} implies easily that
\begin{equation}\label{64}
a q(1) \leq a \int_{0}^{1} q(\xi)d\xi \leq a A.
\end{equation}
By \eqref{33}, we have
\begin{equation}\label{65}
q'(s) = a \int_{0}^{s}q(\xi)d\xi - a sq(s)  - \mu s - F(s) + sF'(s),
\end{equation}
i.e.
\begin{equation}\notag
sF'(s) - F(s) = q'(s) + \mu s + a sq(s) - a \int_{0}^{s}q(\xi)d\xi.
\end{equation}
Thus, for $s \neq 0$ we have
\begin{equation}\label{66}
\left( \frac{F(s)}{s}\right)' = \frac{q'(s)}{s^{2}} + \frac{\mu}{s} + a \left(\frac{1}{s} \int_{0}^{s}q(\xi)d\xi\right)'.
\end{equation}
We pick an $s \geq 1$ and integrate both sides of the equation \eqref{66} from $1$ to $s$.
This results to
\begin{equation}\notag
\frac{F(s)}{s} - F(1) = \int_{1}^{s} \frac{q'(\xi)}{\xi^{2}}d\xi + \mu \ln s + a \frac{1}{s} \int_{0}^{s} q(\xi)d\xi - a \int_{0}^{1} q(\xi)d\xi,
\end{equation}
or
\begin{equation}\label{67}
\int_{1}^{s} \frac{q'(\xi)}{\xi^{2}}d\xi = \frac{F(s)}{s} - F(1) - \mu \ln s - \frac{a}{s} \int_{0}^{s} q(\xi)d\xi + a \int_{0}^{1} q(\xi)d\xi.
\end{equation}
Since $q'(s) < 0$, for all $s \in (0, \infty)$,
\begin{equation}\notag
0 \geq \int_{1}^{s} \frac{q'(\xi)}{\xi^{2}}d\xi \geq \int_{1}^{s} q'(\xi) d\xi \geq \int_{0}^{\infty} q'(\xi) d\xi = \lim_{s \rightarrow \infty} q(s) - q(0) = - A
\end{equation}
hence, \eqref{67} gives
\begin{equation}\notag
0 \geq \frac{F(s)}{s} - F(1) - \mu \ln s - \frac{a}{s} \int_{0}^{s} q(\xi)d\xi + a \int_{0}^{1} q(\xi)d\xi \geq - A,
\end{equation}
or
\begin{equation}\label{68}
F(1) + \mu \ln s + \frac{a}{s} \int_{0}^{s} q(\xi)d\xi - a \int_{0}^{1} q(\xi)d\xi \geq \frac{F(s)}{s} \geq  F(1) + \mu \ln s + \frac{a}{s} \int_{0}^{s} q(\xi)d\xi - a \int_{0}^{1} q(\xi)d\xi - A.
\end{equation}
By \eqref{65} we have
\begin{equation}\label{69}
\frac{F(s)}{s} = \frac{a}{s} \int_{0}^{s}q(\xi)d\xi - a q(s) - \mu  + F'(s) - \frac{q'(s)}{s}.
\end{equation}
Then, \eqref{68} combined with \eqref{69} implies
\begin{equation}\notag
F(1) + \mu \ln s - a \int_{0}^{1} q(\xi)d\xi \geq F'(s) - a q(s) - \mu - \frac{q'(s)}{s}  \geq  F(1) + \mu \ln s - a \int_{0}^{1} q(\xi)d\xi - A,
\end{equation}
or
\begin{equation}\label{70}
\frac{q'(s)}{s} + F(1) + \mu \ln s - a \int_{0}^{1} q(\xi)d\xi \geq F'(s) - a q(s) - \mu \geq \frac{q'(s)}{s} + F(1) + \mu \ln s - a \int_{0}^{1} q(\xi)d\xi - A.
\end{equation}
By \eqref{48} and the fact that $\left\|q' \right\|_{\infty} \geq - q'(s)$ for all $s \in (0, \infty)$, we have
\begin{equation}\notag
- q'(s) \leq \left\|q' \right\|_{\infty} \leq \mu \sqrt{\frac{3 A}{1 + 3 a A}}
\qquad \text{for all }\; s > 0,
\end{equation}
which implies
\begin{equation}\label{72}
- \frac{\mu}{s}\sqrt{\frac{3 A}{1 + 3 a A}} \leq \frac{q'(s)}{s} < 0 < - \frac{q'(s)}{s} \leq \frac{\mu}{s}\sqrt{\frac{3 A}{1 + 3 a A}}.
\end{equation}
Hence, by using \eqref{72} in \eqref{70} we obtain
\begin{align}
\frac{\mu}{s}\sqrt{\frac{3 A}{1 + 3 a A}} + F(1) + \ln s^{\mu} - a \int_{0}^{1} q(\xi)d\xi & > F'(s) - a q(s) - \mu 
\nonumber
\\
&\geq - \frac{\mu}{s}\sqrt{\frac{3 A}{1 + 3 a A}} + F(1) + \ln s^{\mu} - a \int_{0}^{1} q(\xi)d\xi - A
\nonumber
\end{align}
for all $s > 0$. Now, by invoking \eqref{64} the inequalities above give
\begin{equation}\notag
\frac{\mu}{s}\sqrt{\frac{3 A}{1 + 3 a A}} + F(1) + \ln s^{\mu} - a q(1) > F'(s) - a q(s) - \mu \geq - \frac{\mu}{s}\sqrt{\frac{3 A}{1 + 3 a A}} + F(1) + \ln s^{\mu} - a A - A.
\end{equation}
Using the definition of $F(s)$ the above inequalities can be written in the form
\begin{equation}\notag
\frac{\mu}{s}\sqrt{\frac{3 A}{1 + 3 a A}} + F(1) + \ln s^{\mu} - a q(1) > - \frac{1}{3} \ln q(s) - a q(s) - \mu \geq - \frac{\mu}{s}\sqrt{\frac{3 A}{1 + 3 a A}} + F(1) + \ln s^{\mu} - a A - A.
\end{equation}
Recalling \eqref{63}, the above inequalities imply
\begin{align}
3 a q(1) + \ln q(1) &- \frac{3\mu}{s}\sqrt{\frac{3 A}{1 + 3 a A}} + \ln s^{-3 \mu} - 3\mu 
< \ln q(s) + 3a q(s)
\nonumber
\\
&\leq 3 a A + 3A + \frac{3\mu}{s}\sqrt{\frac{3 A}{1 + 3 a A}} + \ln A + \ln s^{-3 \mu} - 3\mu,
\nonumber
\end{align}
which by exponentiation yields \eqref{61}.
\hfill $\blacksquare$

\smallskip

\textbf{Corollary 2.2.}
Let $q(s)$ satisfy \eqref{24}-\eqref{25} (in particular $q(0) = A$).
Then, as a function of A, the quantity
\begin{equation}\label{73}
I(A) := \int_{- \infty}^{\infty} q(s) ds
\end{equation}
is continuous in $(0, \infty)$.

\smallskip

\textit{Proof}.
Let $q(s) = q (s ; A)$ be the unique solution of the problem \eqref{24}-\eqref{25}.
By the standard theorem of ordinary differential equations on continuous dependence on the parameters we have that
$q(s ; A)$ is continuous in A for all $A > 0$.
For fixed $A_{1}, A_{2}$ with $0 < A_{1} < A_{2} < \infty$, the second inequality in \eqref{61}, the monotonicity
of $q$ and the condition $\mu > 1/3$ imply that the family $\left\{q(\cdot \, ; A) : A \in [A_{1}, A_{2}]\right\}$
is dominated by the integrable function $H(s) = h(\left| s \right|)$, $s\in \mathbb{R}$, where
\begin{equation}\notag
h(s) :=
\begin{cases}
A_{2} e^{3 a A_{2}}, \qquad 0 \leq s \leq 1;
\\
\frac{A_{2} e^{3 A_{2} (1 + a)}}{s^{3 \mu} \exp \left[3 \mu  - 3 \mu \sqrt{\frac{3 A_{2}}{1 + 3 a A_{1}}}\right]}, \qquad s \geq 1.
\end{cases}
\end{equation}
Hence, the continuity of $I(A)$ follows by invoking the Dominated Convergence Theorem.
\hfill $\blacksquare$

\smallskip

We are now ready for our main result.

\smallskip

\textbf{Theorem 2.1.} Let $\gamma = -2 / 3$. Then, for
each number $\beta \in (0, \infty)$ there is a self-similar solution of \eqref{7} (which is equivalent to \eqref{A3}), namely a solution $u$
of the form $u(t,x) = t^{-1/3}g(x t^{-1/3})$, where $g(s)$ satisfies \eqref{20}, \eqref{21}, \eqref{22}, and \eqref{23}, such that
\begin{equation}\notag
\beta = K[g] + \frac{a}{2} \Lambda [g].
\end{equation}

\smallskip

\textit{Proof}.
Let $q(s) = q(s ; A)$ be the unique solution of the problem \eqref{24}-\eqref{25} with $\mu = \beta + (1/3)$, that is
\begin{equation}\notag
q''(s)q(s) + a s q'(s)q(s) + \beta q(s) + \frac{1}{3}q(s) + \frac{1}{3}sq'(s) = 0,
\qquad
s \in \mathbb{R},
\end{equation}
\begin{equation}\notag
q(0) = A > 0, \qquad q'(0)=0,
\end{equation}
and set
\begin{equation}\label{74}
Q(A) := \int_{-\infty}^{\infty} q'(s ; A)^2 ds + \frac{a}{2}\int_{-\infty}^{\infty} q(s ; A)^2 ds.
\end{equation}
Then by \eqref{31} of Lemma 2.2
\begin{equation}\notag
Q(A) = \beta \int_{-\infty}^{\infty}q(s ; A) ds = \beta I(A)
\end{equation}
(recall \eqref{73}), hence Corollary 2.2 tells us that $Q(A)$ is continuous on $(0, \infty)$.
Furthermore, by \eqref{54} of Corollary 2.1 we have
\begin{equation}\label{75}
\lim_{A \rightarrow 0^{+}}Q(A) = 0
\qquad \text{and} \qquad
\lim_{A \rightarrow \infty}Q(A) = \infty.
\end{equation}
Thus, $Q(A)$ takes every value between $0$ and $\infty$.
In particular, for each number $\beta \in (0, \infty)$ there is an $A = A_{\beta}$ such that
\begin{equation}\notag
Q(A_{\beta}) = \beta.
\end{equation}
Set $g(s) = q(s ; A_{\beta})$.
Then
\begin{equation}\notag
 K[g] + \frac{a}{2} \Lambda [g] = \int_{-\infty}^{\infty} q'(s ; A)^2 ds + \frac{a}{2}\int_{-\infty}^{\infty} q(s ; A)^2 ds
 = Q(A_{\beta}) = \beta ,
 \end{equation}
hence $g(s)$ satisfies \eqref{21}-\eqref{22}.
Furthermore, by \eqref{31} of Lemma 2.2
\begin{equation}\notag
\int_{-\infty}^{\infty}g(s)ds = \int_{-\infty}^{\infty}q(s ; A_{\beta})ds
= \frac{1}{\beta} \left[\int_{-\infty}^{\infty} q'(s ; A_{\beta})^2 ds + \frac{a}{2}\int_{-\infty}^{\infty} q(s ; A_{\beta})^2 ds \right]
= \frac{1}{\beta}Q(A_{\beta}) = 1,
\end{equation}
and, therefore, $g(s)$ also satisfies \eqref{20}.
\hfill $\blacksquare$

\smallskip

Clearly, all these self-similar solutions $u(t, x)$ are probability density functions on $\mathbb{R}$. A peculiar feature of these solutions is that
they all approach the Dirac delta function $\delta(x)$ as $t \to 0^+$.


\begin{thebibliography}{9}

\bibitem{Bo} Bomze, I., Dynamical aspects of evolutionary stability,
\textit{Monaish. Mathematik} {\bf 110} (1990), 189--206.

%\bibitem{C} Chung, K.L., {\it Lectures from Markov Processes to Brownian Motion},
%A Series of Comprehensive Studies in Mathematics,
%Springer-Verlag New York, Inc., 1982.

\bibitem{C-L} Coddington, E.A. and Levinson, N., {\it Theory
of Ordinary Differential Equations}, Robert E. Krieger Publishing
Company, Malabar, Florida, 1987.

%\bibitem{C-H} Courant, R. and Hilbert, D., {\it Methods of Mathematical
%Physics, Volume I}, Interscience Pub. New York, 1953.
%
%\bibitem{D-S} Dunford, N. and Schwartz, J.T., {\it Linear Operators, Part II:
%Spectral Theory: Self Adjoint Operators in Hilbert Space}, Wiley Interscience, 1963.
%
%\bibitem{D} Durrett, R., \textit{Brownian Motion and Martingales in Analysis},
%Wadsworth, Inc., Belmont, CA, 1984.

\bibitem{H-S} Hofbauer, J. and Sigmund, K., Evolutionary Game Dynamics,
\textit{Bulletin (New Series) of the American Mathematical Society} {\bf 40}, no. 4 (2003),
479--519.

\bibitem{I} Imhof, L.A., The long-run behavior of the stochastic
replicator dynamics, \textit{Ann. Appl. Probab.} {\bf 15}, no. 1B (2005),
1019--1045.

%\bibitem{K-S} Karatzas, I. and Shreve, S.E., \textit{Brownian Motion and Stochastic Calculus}, Second
%Edition, Springer, New York, 1991.

\bibitem{K-P} Kravvaritis, D., Papanicolaou, V.G., and
Yannacopoulos, A.,  Similarity Solutions for a Replicator Dynamics
Equation, \textit{Indiana Univ. Math. Journal} {\bf 57}, no. 4 (2008), 1927--1943.

\bibitem{Oe} Oechssler, J. and Riedel, F., Evolutionary dynamics on infinite strategy spaces,
\textit{Economic Theory}  {\bf 17} (2001), 141--162.

\bibitem{Oe1} Oechssler, J. and Riedel, F., On the dynamic foundation of evolutionary stability in
continuous models, \textit{Journal of Economic Theory} {\bf 107} (2002),
223--252.

\bibitem{P-S} Papanicolaou, V.G. and Smyrlis G.,
Similarity Solutions for a Multidimensional Replicator Dynamics Equation,
\textit{Nonlinear Analysis} {\bf 71} (2009), 3185--3196.

%\bibitem{P} Pazy, A., {\it Semigroups of Linear Operators and applications to Partial Differential Equations},
%Springer-Verlag New York, Inc., 1983.
%
%\bibitem{R-S} Reed, M. and Simon, B., \textit{Methods of Moderm Mathematical Physics I: Functional Analysis},
%Revised and Enlarged Edition, Academic Press, Inc., San Diego, CA, 1980.

%\bibitem{R-Y} Revuz, D. and Yor, M., \textit{Continuous Martingales and Brownian Motion},
%Third Edition, A Series of Comprehensive Studies in Mathematics 293, Springer, New York, 2005.

\bibitem{S} Smith, J. Maynard, {\it Evolution and the Theory of Games},
Cambridge University Press, Cambridge, UK, 1982.

\bibitem{T-J} Taylor, P.D. and Jonker, L.B., Evolutionary Stable
Strategies and Game Dynamics, \textit{Mathematical Biosciences} {\bf 40}
(1978), 145--156.

%\bibitem{W} Weidmann, J., {\it Spectral Theory of Ordinary Differential Operators},
%Lecture Notes in Mathematics, vol. 1258, Springer, 1987.

%\bibitem{Z} Zettl, A., {\it Sturm-Liouville Theory}, Mathematical Surveys and Monographs,
%vol. 121, American Mathematical Society, Providence, RI, 2005.
%
%\bibitem {Ok} {\O}ksendal B., \textit{Stochastic differential Equations}, Springer-Verlag, 1995
%

\end{thebibliography}
\end{document}